\documentclass[12pt]{article}
\usepackage{amscd,amsfonts,amssymb,amsmath,latexsym,array,hhline}
\usepackage[dvips]{graphics}
\makeatletter
\@addtoreset{equation}{section}
\renewcommand{\@begintheorem}[2]{\begin{trivlist}\it
\item[\hspace{\labelsep}{\bf #1\ #2.}]}
\renewcommand{\@opargbegintheorem}[3]{\begin{trivlist}\it
\item[\hspace{\labelsep}{\bf #1\ #2\ (#3).}]}
\renewcommand{\@endtheorem}{\end{trivlist}}
\renewcommand{\@cite}[2]{[{#1\if@tempswa ; #2\fi}]}
\makeatother

\newcommand{\paragr}{\hspace{6mm}}
\newtheorem{Theorem}{\paragr Theorem}[section]
\newtheorem{Lemma}[Theorem]{\paragr Lemma}

\newtheorem{Definition}[Theorem]{\paragr Definition}

\newtheorem{Example}[Theorem]{\paragr Example}

\newcommand{\Proof}{\texttt{Proof}. }
\newcommand{\Remark}{\texttt{Remark}. }
\newcommand{\Iremark}{\texttt{Important remark}. }

\newcommand{\Remarks}{\texttt{Remarks}. }

\newcommand{\End}{~\hfill $\Box\!$}

\newcommand{\skm}{\bigskip}
\newcommand{\skb}{\bigskip}

\newcommand{\al}{\alpha}
\newcommand{\be}{\beta}

\newcommand{\de}{\delta}

\newcommand{\la}{\lambda}

\newcommand{\eps}{\varepsilon}
\renewcommand{\phi}{\varphi}
\renewcommand{\kappa}{\varkappa}
\newcommand{\N}{\mathbb{N}}

\newcommand{\R}{\mathbb{R}}

\newcommand{\B}{\mathcal{B}}
\newcommand{\F}{\mathcal{F}}

\newcommand{\AAA}{\mathcal{A}}

\newcommand{\DDD}{\mathcal{D}}

\newcommand{\RRR}{\mathcal{R}}
\newcommand{\PPP}{\mathcal{P}}
\newcommand{\VVV}{\mathcal{V}}
\newcommand{\EX}{\mathcal{X}}

\newcommand{\EE}{\mathsf{E}}
\newcommand{\PP}{\mathsf{P}}
\newcommand{\QQ}{\mathsf{Q}}

\newcommand{\emp}{\emptyset}

\newcommand{\wt}{\widetilde}
\newcommand{\wl}{\overline}
\newcommand{\ul}{\underline}
\newcommand{\xra}{\xrightarrow}
\newcommand{\da}{\downarrow}

\newcommand{\cond}{\hspace{0.3mm}|\hspace{0.3mm}}

\newcommand{\Law}{\mathop{\rm Law}\nolimits}
\newcommand{\conv}{\mathop{\rm conv}\nolimits}

\newcommand{\argmin}{\mathop{\rm argmin}}

\renewcommand{\inf}{\mathop{\rm inf\rule[-0.8mm]{0mm}{1mm}}}

\newcommand{\Lea}{\Longleftrightarrow}

\newenvironment{mitemize}%
{\begin{list}{$\bullet$}{
\leftmargin=32pt
\rightmargin=0pt
\labelsep=5pt
\labelwidth=20pt
\itemindent=0pt
\topsep=5pt plus 2pt minus 4pt
\partopsep=2pt plus 1pt minus 1pt
\parsep=0pt
\itemsep=0pt}}%
{\end{list}}

\begin{document}
\vspace*{5mm}
\begin{center}\bf
PRICING AND HEDGING

\vspace{1mm}
IN INCOMPLETE MARKETS

\vspace{1mm}
WITH COHERENT RISK
\end{center}

\begin{center}\itshape\bfseries
Alexander S.~Cherny$^*$,\quad Dilip B.~Madan$^{**}$
\end{center}

\begin{center}
\textit{$^*$Moscow State University}\\
\textit{Faculty of Mechanics and Mathematics}\\
\textit{Department of Probability Theory}\\
\textit{119992 Moscow, Russia}\\
\texttt{E-mail: cherny@mech.math.msu.su}\\
\texttt{Webpage: http://mech.math.msu.su/\~{}cherny}
\end{center}

\begin{center}
\textit{$^{**}$Robert H.~Smith School of Business}\\
\textit{Van Munching Hall}\\
\textit{University of Maryland}\\
\textit{College Park, MD 20742}\\
\texttt{E-mail: dmadan@rhsmith.umd.edu}\\
\texttt{Webpage: http://www.rhsmith.umd.edu/faculty/dmadan}
\end{center}

\begin{abstract}
\textbf{Abstract.}
We propose a pricing technique based on coherent risk
measures, which enables one to get finer price intervals
than in the No Good Deals pricing.
The main idea consists in splitting a liability into
several parts and selling these parts to different
agents.
The technique is closely connected with the
convolution of coherent risk measures and
equilibrium considerations.

Furthermore, we propose a way to apply the above
technique to the coherent estimation of the Greeks.

\medskip
\textbf{Key words and phrases.}
CDO,
coherent risk measure,
convolution-based pricing,
extreme measure,
incomplete markets,
factor risk,
maximum-based pricing,
No Strictly Acceptable Opportunities,
risk contribution,
sensitivity coefficients,
valuation measure,
Weighted V@R.
\end{abstract}

\section*{Contents}

\noindent\hbox to \textwidth{\ref{I}\ \ Introduction
\dotfill \pageref{I}}

\noindent\hbox to \textwidth{\ref{CRM}\ \ Basic
Definitions \dotfill \pageref{CRM}}

\noindent\hbox to \textwidth{\ref{BO}\ \ Basic
Operations \dotfill \pageref{BO}}

\noindent\hbox to \textwidth{\ref{MBP}\ \
Maximum-Based Pricing \dotfill \pageref{MBP}}

\noindent\hbox to \textwidth{\ref{CBP}\ \
Convolution-Based Pricing \dotfill \pageref{CBP}}

\noindent\hbox to \textwidth{\ref{SENS}\ \
Sensitivity Coefficients \dotfill \pageref{SENS}}

\noindent\hbox to \textwidth{\ref{SC}\ \
Summary and Conclusion \dotfill \pageref{SC}}

\noindent\hbox to \textwidth{References
\dotfill \pageref{R}}

\newpage
\section{Introduction}
\label{I}

\textbf{1. Pricing and hedging in incomplete markets.}
One of the main problems of the modern finance is:
\textit{How to price derivatives in incomplete markets?}

The arbitrage theory typically provides price
intervals that are unacceptably large.
In order to narrow these intervals, one should use ideas
beyond No Arbitrage.

One of possible methods has been proposed by Carr,
Geman, and Madan~\cite{CGM01}.
Their idea is as follows.
Suppose that we have~$N$ groups of agents and the $n$-th
group assesses the quality of any possible trade by
taking the expectation of its P\&L~$X$ with respect to
\textit{valuation measures}~$\QQ^{nk}$, $k=1,\dots,K^n$.
Thus, a trade is profitable for this group if
$\EE_{\QQ^{nk}}X>0$ for any~$k$.
Natural examples of $\QQ^{nk}$ are given in Section~\ref{BO}.
A trade is called \textit{strictly acceptable} if its
P\&L satisfies $\EE_{\QQ^{nk}}X>0$ for any $n,k$.
The No Strictly Acceptable Opportunities assumption
says that such trades do not exist.
(This is a natural strengthening of the No Arbitrage assumption.)
The fundamental theorem of asset pricing provided in~\cite{CGM01}
states that strictly acceptable opportunities do not
exist if and only if
$$
\bigl(\conv\limits_{n,k}\QQ^{nk}\bigr)\cap\RRR\ne\emp,
$$
where $\RRR$ is the set of risk-neutral measures.
Using this ideology, one can define the upper and lower prices
of a contingent claim~$F$ as
\begin{align*}
\wl V(F)&=\inf\bigl\{x:\exists X\in A\text{ such that }
\inf\limits_{n,k}\EE_{\QQ^{nk}}(X-F+x)>0\bigr\},\\
\ul V(F)&=\sup\bigl\{x:\exists X\in A\text{ such that }
\inf\limits_{n,k}\EE_{\QQ^{nk}}(X+F-x)>0\bigr\},
\end{align*}
where $A$ is the set of P\&Ls that can be obtained from
various trades available on the market.
It follows from the above result that
\begin{align*}
\wl V(F)&=\sup\bigl\{\EE_\QQ F:\QQ\in\bigl(
\conv\limits_{n,k}\QQ^{nk}\bigr)\cap\RRR\bigr\},\\
\ul V(F)&=\inf\bigl\{\EE_\QQ F:\QQ\in\bigl(
\conv\limits_{n,k}\QQ^{nk}\bigr)\cap\RRR\bigr\}.
\end{align*}
Let us remark that this pricing technique is closely
connected with the \textit{No Good Deals} pricing;
see~\cite{BL00}, \cite{CGM04}, \cite{CH01},
\cite{C061}, \cite{CS00}, \cite{JK01}.

The basic idea of this paper is as follows.
Let us call a trade producing a P\&L~$X$
\textit{strictly acceptable for the $n$-th group}
if $\EE_{\QQ^{nk}}X>0$ for any~$k=1,\dots,K^n$.
Let us now call a trade~$X$ \textit{strictly acceptable}
if~$X$ can be represented as $X^1+\dots+X^N$, where
$X^n$ is strictly acceptable for the $n$-th group.
The fundamental theorem of asset pricing provided in this
paper states that there are no strictly acceptable
opportunities if and only if
$$
\Bigl(\bigcap_{n=1}^N\conv\limits_{k=1,\dots,K^n}
\QQ^{nk}\Bigr)\cap\RRR\ne\emp.
$$
Using this ideology, we can define the upper and lower prices
of a contingent claim~$F$ as
\begin{align*}
\wl V(F)
&=\inf\bigl\{x:\exists X\in A,\,Y^1,\dots,Y^N
\text{ such that }\textstyle\sum_n\!Y^n=X-F+x\\
&\hspace*{11mm}\text{ and }\inf_{k}\EE_{\QQ^{nk}} Y^n>0
\text{ for any }n\bigr\},\\
\ul V(F)
&=\sup\bigl\{x:\exists X\in A,\,Y^1,\dots,Y^N
\text{ such that }\textstyle\sum_n\!Y^n=X+F-x\\
&\hspace*{11mm}\text{ and }\inf_{k}\EE_{\QQ^{nk}} Y^n>0
\text{ for any }n\bigr\}.
\end{align*}
In other words, superreplicating a liability consists
of two steps:
\begin{mitemize}
\item[\bf 1.] Trading in the market.
\item[\bf 2.] Splitting a liability into several parts
and selling these parts to different groups at the price~0.
\end{mitemize}
As shown in the paper,
\begin{align*}
\wl V(F)&=\sup\Bigl\{\EE_\QQ F:\QQ\in\Bigl(
\bigcap_{n=1}^N\conv\limits_{k=1,\dots,K^n}\QQ^{nk}\Bigr)
\cap\RRR\Bigr\},\\
\ul V(F)&=\inf\Bigl\{\EE_\QQ F:\QQ\in\Bigl(
\bigcap_{n=1}^N\conv\limits_{k=1,\dots,K^n}\QQ^{nk}\Bigr)
\cap\RRR\Bigr\}.
\end{align*}
Clearly, this technique leads to finer price intervals
than the technique described above.
A remarkable property is: the more
groups are taken into account, the smaller are the fair
price intervals. Note that the technique described above
has exactly the opposite effect.

Let us remark that the proposed pricing technique is closely
connected with the risk sharing (or equilibrium)
problem considered in~\cite{BE05}, \cite{BE06},
\cite{C062}, \cite{JST05}.

\skb
\textbf{2. Sensitivity coefficients.}
Two main goals of asset pricing are:
\begin{mitemize}
\item[\bf 1.] finding fair prices of OTC derivatives;
\item[\bf 2.] finding sensitivity coefficients, which is needed
for risk management.
\end{mitemize}
In the above discussion, we have concentrated on the
first problem. Let us now discuss the second one.

Traditionally, in complete models the fair price of
a contingent claim is a single number, so that the
sensitivity coefficient of the price to the value of an
underlying asset can be defined simply by taking
the partial derivative.
In incomplete markets, we typically have a whole interval
of fair prices. The above efforts were aimed at narrowing
this interval, but typically it is not a one-point set.
Thus, the partial derivative cannot be taken and
additional ideas should be used to define sensitivity
coefficients.
One way to do this was proposed in~\cite[Subsect.~3.3]{CS00}.
In this paper, we take another path and adjust the pricing
technique described above to obtain an interval for a
sensitivity.

\skb
\textbf{3. Coherent risk.}
The pricing technique described above is closely
connected with the notion of a \textit{coherent risk
measure}. This notion was introduced by Artzner,
Delbaen, Eber, and Heath~\cite{ADEH97}, \cite{ADEH99}
as a substitute for V@R, which is notorious for its
drawbacks.
This notion turned out to be extremely useful
not only for the risk measurement purposes, but also for
other financial problems like pricing and optimization
(see~\cite{C061}, \cite{C062} and references therein).
A coherent risk measure is a function on random variables
defined as
\begin{equation}
\label{i1}
\rho(X)=-\inf_{\QQ\in\DDD}\EE_\QQ X,
\end{equation}
where $\DDD$ is a set of probability measures
($X$ means the P\&L of some transaction).
This formula has a clear interpretation: we have a
family~$\DDD$ of probabilistic scenarios;
we calculate the average P\&L under each scenario;
then we take the worst case.

The relationship between the above considerations and
coherent risk measures is obvious.
Define the coherent risk measure of the $n$-th group as
$$
\rho^n(X)=-\inf_{k=1,\dots,K^n}\EE_{\QQ^{nk}}X.
$$
The first pricing technique says that $X$ is strictly
acceptable if $\max_n\rho^n(X)<0$.
For this reason, we call here this technique the
\textit{maximum-based pricing}.
Note that $\max_n\rho^n$ is again a coherent risk
measure.

The second technique says that $X$ is strictly
acceptable if $\rho(X):=\min\sum_n\rho^n(X^n)<0$, where
the minimum is taken over all $X^1,\dots,X^N$ such that
$\sum_n X^n=X$.
The function~$\rho$ is again a coherent risk measure
termed the \textit{convolution of $\rho^1,\dots,\rho^N$}.
For this reason, we call this pricing technique the
\textit{convolution-based pricing}.

Both pricing techniques described above are taking the path
$$
\QQ^{nk}\;\longrightarrow\;\rho^n.
$$
However, one of the most natural ways of constructing
the valuation measures $\QQ^{nk}$ is inverting this
arrow, i.e. we start from a coherent risk measure~$\rho^n$
of the $n$-th group, find the set $\DDD^n$ such that
$\rho^n(X)=-\inf_{\QQ\in\DDD^n}X$, and take~$\DDD^n$
as the set of valuation measures for the $n$-th group.
This is meaningful because typically coherent risk
measures are defined not through~\eqref{i1}, but in a
more direct way.
For example, the best (in our opinion) one-parameter
family of coherent risk measures termed
\textit{Alpha V@R} is defined as:
$$
\rho_\al(X)=-\EE\min_{i=1,\dots,\al} X_i,
$$
where $\al$ is a natural number and $X_1,\dots,X_\al$ are
independent copies of~$X$.
This class of risk measures was introduced in~\cite{CM061}.
As shown in that paper, it is very convenient for the
risk measurement purposes.

The sets of valuation measures might be constructed
not only from the ``pure'' risk measures like Alpha V@R, but
also through some ``derivative'' risk measures like
factor risks or risk contributions.
They might also be obtained through some basic operations
on coherent risks.

In typical cases, the sets
$\DDD^n$ thus obtained are infinite.
Therefore, in this paper we are dealing from the outset
with general sets~$\DDD^n$ rather than
with valuation measures~$\QQ^{nk}$.
Moreover, as opposed to the majority of papers on coherent
risk, we are dealing with unbounded random variables
(clearly, this is necessary for financial applications)
using the framework introduced in~\cite{C061}, \cite{C062}.

\skb
\textbf{4. Structure of the paper.}
In Section~\ref{CRM}, we recall basic definitions
related to coherent risk measures and give several
basic examples.

Section~\ref{BO} has two main goals.
First, we provide natural examples of valuation
measures.
These are based on various transformations of coherent
risk measures like factor risks or risk contributions.
Second, we describe two basic operations on coherent
risks: maximum and convolution.
The first operation is at the basis of the pricing
technique of Section~\ref{MBP};
the second one is at the basis of the pricing technique
of Section~\ref{CBP}.

In Section~\ref{MBP}, we consider the maximum-based
pricing. This is, in fact, the pricing technique
proposed in~\cite{CGM01} and extended to a more
general framework in~\cite{C061}.
This section is included just for the reader's
convenience.

The main results of this paper are presented in
Section~\ref{CBP}. For the convolution-based
pricing technique, we establish the fundamental
theorem of asset pricing and find the form of the
price intervals.
On the mathematical side, our statements easily
follow from the results of~\cite{C061}, \cite{C062}.
We illustrate this technique with an example of
calculating the prices and the superreplication strategy
for the case, where different groups are using coherent
risk measures from the class Weighted V@R.
The quantitative effect arising in this example is very
similar to how CDOs are arranged.
We also discuss the empirical estimation of fair price
intervals provided by other natural classes of valuation measures.
Furthermore, we study the liquidity effects following
the ideas of~\cite{CGM01} and the technique of~\cite{C062}.

In Section~\ref{SENS}, we adjust the convolution-based
pricing technique to the assessment of sensitivity
coefficients.

\section{Basic Definitions}
\label{CRM}

\textbf{1. Coherent risk measures.}
Let $(\Omega,\F,\PP)$ be a probability space.
Recall that $L^\infty$ is the space of bounded random
variables on $(\Omega,\F,\PP)$.
The following definition was introduced in~\cite{ADEH97},
\cite{ADEH99}, \cite{D02}.

\begin{Definition}\rm
\label{CRM1}
A \textit{coherent risk measure on $L^\infty$} is a
map $\rho\colon L^\infty\to\R$ satisfying the properties:
\begin{mitemize}
\item[(a)] (Subadditivity) $\rho(X+Y)\le\rho(X)+\rho(Y)$;
\item[(b)] (Monotonicity) If $X\le Y$, then
$\rho(X)\ge\rho(Y)$;
\item[(c)] (Positive homogeneity) $\rho(\la X)=\la\rho(X)$ for
$\la\in\R_+$;
\item[(d)] (Translation invariance)
$\rho(X+m)=\rho(X)-m$ for $m\in\R$;
\item[(e)] (Fatou property) If $|X_n|\le1$,
$X_n\xra{\PP}X$, then $\rho(X)\le\liminf_n\rho(X_n)$.
\end{mitemize}
\end{Definition}

The theorem below is the basic representation theorem.
It was established in~\cite{ADEH99} for
the case of a finite $\Omega$ (in this case the axiom~(e)
is not needed) and in~\cite{D02} for the general case.
By~$\PPP$ we will denote the set of probability measures
that are absolutely continuous with respect to~$\PP$.
Throughout the paper, we identify measures from~$\PPP$
(these are typically denoted by~$\QQ$)
with their densities with respect to~$\PP$
(these are typically denoted by~$Z$).

\begin{Theorem}
\label{CRM2}
A function~$\rho$ satisfies
conditions {\rm(a)--(e)} if and only if there exists a
non-empty set $\DDD\subseteq\PPP$ such that
$$
\rho(X)=-\inf_{\QQ\in\DDD}\EE_\QQ X,\quad X\in L^\infty.
$$
\end{Theorem}

\Iremark
Suppose that an agent is assessing risk with a coherent
risk measure~$\rho(X)=-\inf_{\QQ\in\DDD}\EE_\QQ X$.
The corresponding \textit{risk-adjusted performance}
is naturally defined as $p(X):=\EE X-\la\rho(X)$, where
$\la$ is a positive parameter. Then
$$
(1+\la)^{-1}p(X)=\inf_{\QQ\in\DDD'}\EE_\QQ X,
$$
where
$$
\DDD'=\frac{1}{1+\la}\,\PP+\frac{\la}{1+\la}\,\DDD.
$$
The negative of a coherent risk measure is called a
\textit{coherent utility function} (in some respects
this object is more convenient than a coherent risk
measure because it eliminates numerous minus signs).
Thus, $(1+\la)^{-1}p(X)$ is a coherent utility.
This stability property shows that coherent
risk measures can be used not only to measure risk,
but to measure risk-adjusted performance, i.e. utility
as well.\End

\skm
Throughout this paper, we deal with coherent risk
measures on the space $L^0$ of all random variables.
The following definition was introduced in~\cite{C061}.

\begin{Definition}\rm
\label{CRM3}
A \textit{coherent risk measure on $L^0$} is a map
$\rho\colon L^0\to[-\infty,\infty]$ defined as
\begin{equation}
\label{crm1}
\rho(X)=-\inf_{\QQ\in\DDD}\EE_\QQ X,\quad X\in L^0,
\end{equation}
where $\DDD$ is a non-empty subset of $\PPP$ and
$\EE_\QQ X$ is understood as $\EE_\QQ X^+-\EE_\QQ X^-$
($X^+=\max\{X,0\}$, $X^-=\max\{-X,0\}$)
with the convention $\infty-\infty=-\infty$.
\end{Definition}

The set $\DDD$, for which~\eqref{crm1} is true, is not
unique. However, there exists the largest such set.
It is given by $\{\QQ\in\PPP:
\EE_\QQ X\ge-\rho(X)\text{ for any }X\}$.

\begin{Definition}\rm
\label{CRM4}
The largest set, for which~\eqref{crm1} is true, is
called the \textit{determining set} of~$\rho$.
\end{Definition}

\Iremark
Let $\DDD$ be a subset of~$\PPP$.
Define a coherent risk measure~$\rho$ by~\eqref{crm1}.
The determining set of~$\rho$ might be strictly
larger than~$\DDD$. However, if $\DDD$ is convex and
$L^1$-closed, then $\DDD$ is the determining set of~$\rho$.
Indeed, suppose that the determining set~$\wt\DDD$
is larger than~$\DDD$. Choose $\QQ_0\in\wt\DDD\setminus\DDD$.
By the Hahn--Banach theorem, there exists
$X\in L^\infty$ such that
$\EE_{\QQ_0}X<\inf_{\QQ\in\DDD}\EE_\QQ X=-\rho(X)$,
which is a contradiction.\End

\skm
For more information on coherent risk measures, we refer
to~\cite{CM061}, \cite{D05}, and~\cite[Sect.~4]{FS04}.

\skm
\textbf{2. Examples.}
Let us give examples of four most natural classes
of coherent risk measures.

\begin{Example}[Tail V@R]\rm
\label{CRM5}
\textit{Tail V@R of order} $\la\in(0,1]$
(the terms \textit{Average V@R},
\textit{Conditional V@R}, \textit{Expected Shortfall},
and \textit{Expected Tail Loss} are also used)
is the coherent risk measure~$\rho_\la$ corresponding
to the determining set
$$
\DDD_\la=\Bigl\{\QQ\in\PPP:\frac{d\QQ}{d\PP}
\le\la^{-1}\Bigr\}.
$$

If $X$ has a continuous distribution, then
$$
\rho_\la(X)=-\EE(X\cond X\le q_\la(X)),
$$
where $q_\la(X)$ is the $\la$-quantile of~$X$.
This motivates the term Tail V@R.

For a detailed study of this risk measure, we refer
to~\cite{AT02}, \cite[Sect.~2]{CM061},
\cite[Sect.~4.4]{FS04}.\End
\end{Example}

\begin{Example}[Weighted V@R]\rm
\label{CRM6}
Let $\mu$ be a probability measure on $(0,1]$.
\textit{Weighted V@R with the weighting measure~$\mu$}
(the term \textit{spectral risk measure} is also used)
is the coherent risk measure~$\rho_\mu$ defined as
$$
\rho_\mu(X)=\int_{(0,1]}\rho_\la(X)\mu(d\la).
$$
(One can check that this is indeed a coherent risk measure.)

Weighted V@R admits several equivalent representations.
One of the most convenient representations is:
\begin{equation}
\label{crm2}
\rho_\mu(X)=-\int_0^1 q_x(X)\psi_\mu(x)dx,
\end{equation}
where
\begin{equation}
\label{crm3}
\psi_\mu(x)=\int_{[x,1]}\la^{-1}\mu(d\la),
\quad x\in[0,1].
\end{equation}
In particular, suppose that $X$ takes on values
$x_1,\dots,x_T$ with probabilities $p_1,\dots,p_T$.
Let $x_{(1)},\dots,x_{(T)}$ be the numbers
$x_1,\dots,x_T$ in the increasing order and let $n(i)$
be the number such that $x_{(i)}=x_{n(i)}$. Then
\begin{equation}
\label{crm4}
\rho_\mu(X)=-\sum_{t=1}^T x_{n(t)}\int_{z_{t-1}}^{z_t}
\psi_\mu(x)dx,
\end{equation}
where $z_t=\sum_{i=1}^t p_{n(i)}$.
This representation is convenient for the empirical
estimation of~$\rho$.

The determining set $\DDD_\mu$ of $\rho_\mu$ admits the
following representations:
\begin{equation}
\begin{split}
\label{crm5}
\DDD_\mu
&=\{\QQ\in\PPP:\QQ(A)\le\Psi_\mu(\PP(A))\;
\text{\rm for any }A\in\F\}\\
&=\Bigl\{Z\in L^0:Z\ge0,\;\EE_\PP Z=1,\text{ and }
\int_{1-x}^1\!\!q_s(Z)ds\le\Psi_\mu(x)\;\forall
x\in[0,1]\Bigr\}\\
&=\{Z\in L^0:Z\ge0,\;\EE Z=1,\;\text{\rm and }
\EE(Z-x)^+\le\Phi_\mu(x)\;\forall x\in\R_+\},
\end{split}
\end{equation}
where
\begin{align*}
\Psi_\mu(x)&=\int_0^x\psi_\mu(y)dy,\quad x\in[0,1],\\
\Phi_\mu(x)&=\sup_{y\in[0,1]}(\Psi_\mu(y)-xy),\quad x\in\R_+
\end{align*}
(see~\cite[Th.~4.6]{C05e}, \cite[Th.~4.73]{FS04},
\cite[Th.~1.53]{S05}).

For a detailed study of Weighted V@R, we refer
to~\cite{A02}, \cite{A04}, \cite{C05e},
\cite[Sect.~2]{CM061}, \cite[Sect.~4.6, 4.7]{FS04}.\End
\end{Example}

\begin{Example}[Beta V@R]\rm
\label{CRM7}
Let $\al\in(-1,\infty)$, $\be\in(-1,\al)$.
\textit{Beta V@R with parameters $\al,\be$} is the
Weighted V@R with the weighting measure
$$
\mu_{\al,\be}(dx)={\rm B}(\be+1,\al-\be)^{-1}x^\be(1-x)^{\al-\be-1}dx,
\quad x\in[0,1].
$$
As shown in~\cite{CM061}, for $\al,\be\in\N$, Beta V@R
admits the following simple representation
$$
\rho_{\al,\be}(X)=
-\EE\Bigl[\frac{1}{\be}\sum_{i=1}^\be X_{(i)}\Bigr],
$$
where $X_{(1)},\dots,X_{(\al)}$ are the order statistics
obtained from independent copies $X_1,\dots,X_\al$ of~$X$.
This representation provides a very convenient way
for the empirical estimation of~$\rho_{\al,\be}$.

For a detailed study of this risk measure,
see~\cite[Sect.~2]{CM061}.\End
\end{Example}

\begin{Example}[Alpha V@R]\rm
\label{CRM8}
\textit{Alpha V@R} is obtained from Beta V@R by fixing
$\be=1$. Clearly, if $\al\in\N$, then
$$
\rho_\al(X)=-\EE\min_{i=1,\dots,\al}X_i,
$$
where $X_1,\dots,X_\al$ are independent copies of~$X$.\End
\end{Example}

In our opinion, the most important classes of coherent
risk measures are: Alpha V@R, Beta V@R, and Weighted V@R.

\skm
\textbf{3. $L^1$-spaces.}
For technical purposes, we need to recall the
definition of the \textit{strong $L^1$-space} associated
with a coherent risk measure~$\rho$:
$$
L_s^1(\DDD)=\Bigl\{X\in L^0:\lim_{n\to\infty}
\sup_{\QQ\in\DDD}\EE_\QQ|X|I(|X|>n)=0\Bigr\},
$$
where $\DDD$ is the determining set of~$\rho$.

\begin{Example}\rm
\label{CRM9}
\textbf{(i)} For Weighted V@R,
$$
L_s^1(\DDD_\mu)=\{X\in L^0:\rho(X)<\infty,\,\rho(-X)<\infty\}
$$
(see~\cite[Subsect.~2.2]{C061}).
The right-hand side of this equality was called
in~\cite{C061} the \textit{weak $L^1$-space}.
It has a clear financial interpretation: this is the
set of random variables such that their risk is finite
and the risk of their negatives is finite.

\textbf{(ii)} For Beta V@R with $\be>0$
(in particular, for Alpha V@R), $L_s^1=L^1(\PP)$
(see~\cite[Sect.~3]{CM061}).\End
\end{Example}

\section{Basic Operations}
\label{BO}

\textbf{1. Factor risks.}
Suppose that an agent is assessing risk using a coherent
risk measure~$\rho$.
As proposed in~\cite{CM061}, it is important to consider
not only the pure risk $\rho(X)$, but also the
\textit{factor risks}
$$
\rho^f(X;Y^m)
=-\inf_{\QQ\in\EE(\DDD\cond Y^m)}\EE_\QQ X,
\quad m=1,\dots,M.
$$
Here $Y^1,\dots,Y^M$ are the main market factors affecting
risk like the price of oil, the S\&P~500 index, or
the credit spread (to be more precise, $Y^m$ is the
increment of the $m$-th factor over the unit time period),
$X$ means the P\&L produced by some portfolio over this
period, and
$$
\EE(\DDD\cond Y^m):=\{\EE(Z\cond Y):Z\in\DDD\},
$$
where $\DDD$ is the determining set of~$\rho$.
Thus, $\rho^f(\,\cdot\,;Y^m)$ is again a coherent risk
measure with the determining set $\EE(\DDD\cond Y^m)$.
As shown in~\cite{CM061}, under minor technical
assumptions,
\begin{equation}
\label{bo1}
\rho^f(X;Y^m)=\rho(\EE(X\cond Y^m)).
\end{equation}
If an agent is using these measures to assess risk,
he/she might take $\conv_m\EE(\DDD\cond Y^m)$ as the
set of his/her valuation measures.

Another opportunity is that the agent is assessing risk
using one \textit{multi-factor} risk measure
$$
\rho^f(X;Y^1,\dots,Y^M)
:=-\inf_{\QQ\in\EE(\DDD\cond Y^1,\dots,Y^M)}\EE_\QQ X
=\rho(\EE(X\cond Y^1,\dots,Y^M)).
$$
Then he/she might take $\EE(\DDD\cond Y^1,\dots, Y^M)$
as the set of his/her valuation measures.

\skm
\textbf{2. Risk contribution.}
Measurement of pure risk/utility is meaningful only for
a ``poor'' agent, i.e. an agent without a large endowment.
A wealthy agent (for example, a big company) already
has a large portfolio that produces a random P\&L~$W$.
For such an agent, it is reasonable to assess any trade~$X$
as $\rho(W+X)-\rho(W)$ rather than as $\rho(X)$.
In other words, the quantity of interest is the
\textit{risk contribution} of~$X$ to~$W$.
According to the definition introduced in~\cite{C061},
the risk contribution is
$$
\rho^c(X;W)=-\inf_{\QQ\in\EX_\DDD(W)}\EE_\QQ X,
$$
where $\EX_\DDD(W)$ is the set of \textit{extreme
measures} defined as
$$
\EX_\DDD(W)=\bigl\{\QQ\in\DDD:\EE_\QQ W
=\inf_{\QQ\in\DDD}\EE_\QQ W\in(-\infty,\infty)\bigr\}.
$$
The relevance of this definition is justified by the
following relation (see~\cite[Subsect.~2.5]{C061}):
$$
\rho^c(X;W)
=\lim_{\eps\da0}\eps^{-1}(\rho(W+\eps X)-\rho(W)).
$$
Thus, $\EX_\DDD(W)$ is a natural candidate for the set
of valuation measures for a wealthy agent employing
coherent risk.

If an agent is using the classical expected utility
$\EE U(X)$ to assess the quality of his/her position
(here $U\colon\R\to\R$ is a concave increasing function),
then there exists his/her ``personal'' measure,
with which he/she assesses the quality of any possible trade.
This measure is given by $\QQ=cU'(W_1)\PP$, where
$W_1$ is the agent's wealth at the terminal date
and $c$ is the normalizing constant.
The role of this measure is seen from the equality
$$
\lim_{\eps\da0}\eps^{-1}(\EE U(W_1+\eps X)-\EE U(W_1))
=\EE XU'(W_1)
=c^{-1}\EE_\QQ X.
$$
If $X$ is the P\&L produced by some trade and
$X$ is small as compared to~$W$, then $X$ is profitable
for the agent if and only if $\EE_\QQ X>0$.

Thus, the notion of an extreme measure serves as the
coherent counterpart of $cU'(W_1)\PP$.
The set $\EX_\DDD(W)$ is typically a singleton as seen
from the examples below.

\begin{Example}\rm
\label{BO1}
\textbf{(i)} If $W\in L_s^1(\DDD_\mu)$ has a
continuous distribution, then $\EX_{\DDD_\mu}(W)$
consists of a unique measure $\QQ_\mu(W)=\psi_\mu(F(W))\PP$,
where $\psi_\mu$ is given by~\eqref{crm3} and $F$ is the
distribution function of~$W$ (for the proof,
see~\cite[Sect.~6]{C05e}).

\textbf{(ii)} Let $\Omega=\{1,\dots,T\}$ and $W(t)=w_t$,
where all $w_n$ are different.
Let $w_{(1)},\dots,w_{(T)}$ be the numbers
$w_1,\dots,w_T$ in the increasing order
and let $n(i)$ be the number such that $w_{(i)}=w_{n(i)}$.
Then $\EX_{\DDD_\mu}(W)$ consists
of a unique measure $\QQ_\mu(W)$ given by
$$
\QQ_\mu(W)\{n(t)\}=\int_{z_{t-1}}^{z_t}\psi_\mu(x)dx,
$$
where $\psi_\mu$ is given by~\eqref{crm3} and
$z_t=\sum_{i=1}^t\PP\{n(i)\}$ (for the proof,
see~\cite[Sect.~5]{CM061}).\End
\end{Example}

\textbf{3. Factor risk contribution.}
The risk contribution technique can be combined with
the factor risk technique.
Namely, the \textit{factor risk contribution} is
defined as
$$
\rho^{fc}(X;Y;W)
=-\inf_{\QQ\in\EX_{\EE(\DDD\cond Y)}(W)}\EE_\QQ X.
$$
Here $Y$ means the increment of one or several market
factors over the unit time period.
As shown in~\cite{CM061}, under minor technical
assumptions,
\begin{equation}
\label{bo2}
\rho^{fc}(X;Y;W)=\rho^c(\EE(X\cond Y);\EE(W\cond Y)).
\end{equation}
Thus, $\EX_{\EE(\DDD\cond Y)}(W)$ is a natural candidate
for the set of valuation measures for a wealthy agent
employing factor risk.

\begin{Example}\rm
\label{BO2}
Combining representation~\eqref{crm5} with the Jensen
inequality, we see that
\begin{align*}
\EE(\DDD_\mu\cond Y)
&=\Bigl\{\phi(Y):\phi\ge0,\;
\int_{\R^d}\phi(x)\wt\PP(dx)=1,\;\text{\rm and }\\
&\hspace*{7mm}\int_{\R^d}(\phi(x)-K)^+\wt\PP(dx)\le\Phi_\mu(K)\;
\forall K\in\R_+\Bigr\}\\
&=\{\phi(Y):\phi\in\wt\DDD_\mu\},
\end{align*}
where $\wt\PP=\Law Y$ ($Y$ is $d$-dimensional) and
$\wt\DDD_\mu$ is the determining set of Weighted
V@R defined on the space
$(\wt\Omega,\wt\F,\wt\PP)=(\R^d,\B,\wt\PP)$.
Let $W\in L_s^1(\DDD_\mu)$.
Set $g(y)=\EE(W\cond Y=y)$ (recall from Example~\ref{CRM9}
that $L_s^1(\DDD_\mu)\subseteq L^1$). Then
$$
\EE\phi(Y)W
=\EE\phi(Y)g(Y)
=\int_{\R^d}\phi(y)g(y)\wt\PP(dy),\quad\phi\in\wt\DDD_\mu,
$$
so that
$$
\EX_{\EE(\DDD_\mu\cond Y)}(W)
=\{\phi(Y):\phi\in\wt\EX_{\wt\DDD_\mu}(g)\}.
$$

\vspace{-8mm}
\hfill$\Box$
\end{Example}

\vspace{3mm}
For more information on factor risks and risk contributions,
we refer to~\cite{C061}, \cite{CM061}.

\skb
\textbf{4. Maximum of coherent risks.}
Let $\rho^1,\dots,\rho^N$ be coherent risk measures
with the determining sets $\DDD^1,\dots,\DDD^N$.
We assume that each $\DDD^n$ is $L^1$-closed and
uniformly integrable.
This assumption is very natural: as seen from~\eqref{crm5},
it is satisfied by the determining set of Weighted V@R.
If $\DDD$ satisfies this assumption, then
$\EE(\DDD\cond Y)$ also has this property
(for the proof, see~\cite[Sect.~3]{C061}).
It is easy to check that if $\DDD$ satisfies this
assumption and $W\in L_s^1(\DDD)$, then
$\EX_\DDD(W)$ also has this property.

The \textit{maximum} of $\rho^1,\dots,\rho^N$ is
defined as
$$
\rho(X)=-\inf_{\QQ\in\conv_n\DDD^n}\EE_\QQ X.
$$
According to the remark following Definition~\ref{CRM4},
$\conv_n\DDD^n$ is the determining set of~$\rho$.
Obviously,
$$
\rho(X)=\max_{n=1,\dots,N}\rho^n(X),
$$
which motivates the term.

The operation of taking maximum is useful if an agent
has several alternative coherent ways to assess risk
and is on the safe side provided that all these risks
of his/her position are negative.
Another possible interpretation is: we have a group
of agents employing different coherent risk measures;
then the maximum is the most liberal one.

\skb
\textbf{5. Convolution of coherent risks.}
Let $\rho^1,\dots,\rho^N$ be coherent risk measures
with the determining sets $\DDD^1,\dots,\DDD^N$.
We assume that each $\DDD^1$ is closed and
uniformly integrable and $\bigcap_n\!\DDD^n\ne\emp$.

The \textit{convolution} of $\rho^1,\dots,\rho^N$ is
defined as
$$
\rho(X)=-\inf_{\QQ\in\bigcap_n\!\DDD^n}\EE_\QQ X.
$$
According to the remark following Definition~\ref{CRM4},
$\bigcap_n\!\DDD^n$ is the determining set of~$\rho$.
It follows from~\cite[Th.~4.2]{C062} that, for
$X\in L_s^1:=\bigcap_n\!L_s^1(\DDD^n)$,
\begin{equation}
\label{bo3}
\rho(X)=\inf_{X^n\in L_s^1,\,\sum_n\!X^n=X}
\sum_{n=1}^N\rho^n(X^n).
\end{equation}

\begin{Example}\rm
\label{BO3}
\textbf{(i)} It is seen from~\eqref{crm5} that
$$
\bigcap_{n=1}^N\DDD_{\mu^n}
=\bigl\{\QQ\in\PPP:\QQ(A)\le\min_n\Psi_{\mu^n}(\PP(A))\;
\text{\rm for any }A\in\F\bigr\}.
$$
The function $\Psi=\min_n\Psi_{\mu^n}$ is increasing,
concave, continuous, $\Psi(0)=0$, and $\Psi(1)=1$.
Hence, $\Psi=\Psi_\mu$ with $\mu(dx)=-x\Psi''(dx)$,
where $\Psi''$ is the second derivative of~$\Psi$ taken
in the sense of distributions, i.e. it is the measure
on $(0,1]$ defined by $\Psi''((a,b]):=\Psi'_+(b)-\Psi'_+(a)$,
where $\Psi'_+$ is the right-hand derivative.
Thus, the convolution of Weighted V@Rs is again a
Weighted V@R.

\textbf{(ii)} The maximum of Weighted V@Rs need not be
a Weighted V@R as shown by the following example.
Consider $\mu^1=1/2\,\de_{1/3}+1/2\,\de_1$,
$\mu^2=\de_{2/3}$.
We have
$$
\Psi_{\mu^1}(x)=\begin{cases}
2x,&0\le x\le\frac{1}{3},\\
\frac{1}{2}x+\frac{1}{2},&\frac{1}{3}\le x\le 1,
\end{cases}\qquad
\Psi_{\mu^2}(x)=\begin{cases}
\frac{3}{2}x,&0\le x\le\frac{2}{3},\\
1,&\frac{2}{3}\le x\le1.
\end{cases}
$$
Suppose that the maximum of $\rho_{\mu^1},\rho_{\mu^2}$
has the form $\rho_\mu$ with some~$\mu$.
As seen from~\eqref{crm5}, $\Psi_\mu$ should coincide
with the minimal concave majorant of $\Psi_{\mu^1}$ and
$\Psi_{\mu^2}$, i.e.
$$
\Psi_\mu(x)=\begin{cases}
2x,&0\le x\le\frac{1}{3},\\
x+\frac{1}{3},&\frac{1}{3}\le x\le\frac{2}{3},\\
1,&\frac{2}{3}\le x\le1.
\end{cases}
$$
This means that
$\mu=1/3\,\de_{1/3}+2/3\,\de_{2/3}$.
Now, by considering a random variable taking on the values
$-1$, $0$, and $1000$ with probability $1/3$, we see
that $\rho_\mu$ is not the maximum of
$\rho_{\mu^1},\rho_{\mu^2}$.\End
\end{Example}

For more information on operations on coherent risk
measures, we refer to~\cite[Sect.~5]{D05}.

\skb
\textbf{6. Valuation measures.}
Table~1 summarizes various possible ways of constructing
valuation measures.

The first 3 lines correspond to a ``poor'' agent using
one of 3 possible coherent ways of risk measurement.
As indicated in the remark following Theorem~\ref{CRM2},
the negatives of coherent risk measures might be used
to measure utility.

The next 3 lines correspond to a wealthy agent who
measures risk/utility in a coherent way,
while the 7th line corresponds to a wealthy agent
employing expected utility.

Finally, the 8th line corresponds to the case, where
an agent is employing various risk/ulitity measurement
techniques (that yield the sets $\VVV^1,\dots,\VVV^K$
of valuation measures) or we have a group consisting
of several agents, each with his/her own set~$\VVV^k$
of valuation measures.

\begin{figure}
\hspace*{8mm}\begin{tabular}{|p{66mm}|p{30mm}|p{35mm}|}
\hline
\hspace*{21mm}\mbox{\bf Risk/utility}
\hspace*{8.5mm}\mbox{\bf measurement technique}&
\hspace*{8.5mm}\raisebox{-2mm}{\bf Inputs}\rule{0mm}{5mm}&
\hspace*{7.5mm}\mbox{\bf Valuation}
\hspace*{8mm}\mbox{\bf measures}\\
\hhline{|=|=|=|}
\rule{0mm}{5mm}Coherent risk/utility&
$\rho$&
$\DDD$\\[1mm]
\hline
\rule{0mm}{5mm}Factor risk/utility&
$\rho,Y^1,\dots,Y^M$&
$\conv_m\EE(\DDD\cond Y^m)$\\[1mm]
\hline
\rule{0mm}{5mm}Multi-factor risk/utility&
$\rho,Y^1,\dots,Y^M$&
$\EE(\DDD\cond Y^1,\dots,Y^M)$\\[1mm]
\hline
\rule{0mm}{5mm}Risk/utility contribution&
$\rho,W$&
$\EX_\DDD(W)$\\[1mm]
\hline
\rule{0mm}{5mm}Factor risk/utility contribution&
$\rho,Y^1,\dots,Y^M,W$&
$\conv_m\EX_{\EE(\DDD\cond Y^m)}(W)$\\[1mm]
\hline
\rule{0mm}{5mm}Multi-factor risk/utility contribution&
$\rho,Y^1,\dots,Y^M,W$&
$\EX_{\EE(\DDD\cond Y^1,\dots,Y^M)}(W)$\\[1mm]
\hline
\rule{0mm}{5mm}Expected utility&
$U,W_1$&
$cU'(W_1)\PP$\\[1mm]
\hline
\rule{0mm}{5mm}Combination of several techniques&
$\VVV^1,\dots,\VVV^K$&
$\conv_k\VVV^k$\\[1mm]
\hline
\end{tabular}

\vspace{5mm}
\hspace*{41.5mm}{\small\textbf{Table~1.} Construction of valuation
measures}
\end{figure}

\section{Maximum-Based Pricing}
\label{MBP}

Let $(\Omega,\F,\PP)$ be a probability space
and $\VVV^1,\dots,\VVV^N$ be convex $L^1$-closed
uniformly integrable subsets of~$\PPP$.
From the financial point of view, $\VVV^n$ is the set
of valuation measures used by the $n$-th group.
Let $A$ be a convex subset of
$L_s^1:=\bigcap_n\!L_s^1(\VVV^n)$.
From the financial point of view, this is the set of
various discounted P\&Ls that can be obtained by
various trading operations over the unit time period.
For example, if the agents can trade the assets
$1,\dots,d$ and the $i$-th asset produces the discounted
P\&L~$X^i$, then
$$
A=\Bigl\{\sum_{i=1}^d h^iX^i:(h^1,\dots,h^d)\in H\Bigr\},
$$
where $H$ is a portfolio constraint (if there are no
constraints, then $H=\R^d$).

\begin{Definition}\rm
\label{MBP1}
\textbf{(i)} The set of \textit{strictly acceptable
opportunities for the $n$-th group} is
$$
\AAA^n:=\bigl\{X\in L_s^1:\inf_{\QQ\in\VVV^n}\EE_\QQ X>0\bigr\}.
$$

\textbf{(ii)} The set of \textit{strictly acceptable
opportunities} is
$$
\AAA=\{X\in L_s^1:X\in\AAA^n\;\forall n\}.
$$
\end{Definition}

\begin{Lemma}
\label{MBP2}
Let $\rho^n$ be the coherent risk measure with the
determining set~$\VVV^n$ and let $\rho$ be the maximum
of $\rho^1,\dots,\rho^N$. Then
$$
\AAA
=\bigcap_{n=1}^N\AAA^n
=\{X\in L_s^1:\rho(X)<0\}.
$$
\end{Lemma}

This statement is trivial.

\begin{Definition}\rm
\label{MBP3}
The model satisfies the \textit{No Strictly Acceptable
Opportunities} (\textit{NSAO}) condition if $A\cap\AAA=\emp$.
\end{Definition}

\begin{Definition}\rm
\label{MBP4}
A \textit{risk-neutral measure} is a measure $\QQ\in\PPP$
such that $\EE_\QQ X\le0$ for any $X\in A$
(the expectation $\EE_\QQ X$ here is understood as
$\EE_\QQ X^+-\EE_\QQ X^-$ with the convention
$\infty-\infty=-\infty$).

The set of risk-neutral measures will be denoted by~$\RRR$.
\end{Definition}

\begin{Theorem}[FTAP]
\label{MBP5}
The NSAO condition is satisfied if and only if
$(\conv_n\VVV^n)\cap\RRR\ne\emp$.
\end{Theorem}

\Proof
By Lemma~\ref{MBP2}, the NSAO condition coincides with
the \textit{No Good Deals} (\textit{NGD}) condition
defined in~\cite{C061} and applied to the coherent risk
measure~$\rho$ with the determining
set $\VVV=\conv_n\VVV^n$. Now, the result follows
from~\cite[Th.~3.4]{C061}.\End

\begin{Definition}\rm
\label{MBP6}
Let $F\in L^0$ be the discounted
payoff of some contingent claim.
The \textit{upper} and \textit{lower prices} of~$F$ are
defined as
\begin{align*}
\wl V(F)&=\inf\{x\in\R:\exists X\in A\text{ such that }
X-F+x\in\AAA\},\\
\ul V(F)&=\sup\{x\in\R:\exists X\in A\text{ such that }
X+F-x\in\AAA\}.
\end{align*}
\end{Definition}

\begin{Theorem}[Pricing]
\label{MBP7}
If $A$ is a cone and $F\in L_s^1$, then
\begin{align*}
\wl V(F)&=\sup\{\EE_\QQ F:\QQ\in(\conv_n\VVV^n)\cap\RRR\},\\
\ul V(F)&=\inf\{\EE_\QQ F:\QQ\in(\conv_n\VVV^n)\cap\RRR\}.
\end{align*}
\end{Theorem}

\Proof
Take $x_0\in\R$, set
$A(x_0)=A+\{h(x_0-F):h\in\R_+\}$, and denote the corresponding
set of risk-neutral measures by $\RRR(A(x_0))$.
Clearly, the set
$\{x:\exists X\in A\text{ such that }X-F+x\in\AAA\}$
is an open ray. Using Theorem~\ref{MBP5}, we can write
\begin{align*}
\wl V(F)\ge x_0
&\;\Lea\;\not\!\exists X\in A\text{ such that }X-F+x_0\in\AAA\\
&\;\Lea\;A(x_0)\cap\AAA=\emp\\
&\;\Lea\;(\conv_n\VVV^n)\cap\RRR(A(x_0))\ne\emp\\
&\;\Lea\;\exists\QQ\in(\conv_n\VVV^n)\cap\RRR
\text{ such that }\EE_\QQ F\ge x_0.
\end{align*}
This yields the formula for $\wl V(F)$.
The representation of $\ul V(F)$ is proved similarly.\End

\skm
\Remarks
(i) The above theorem is formally true if the NSAO is violated.
In this case $\wl V(F)=-\infty$ and $\ul V(F)=\infty$.

(ii) The above argument shows that there exist
$\wl\QQ,\ul\QQ\in(\conv_n\VVV^n)\cap\RRR$ such that
$\EE_{\wl\QQ}F=\wl V(F)$, $\EE_{\ul\QQ}(F)=\ul V(F)$.
This is in contrast with the No Arbitrage technique.\End

\section{Convolution-Based Pricing}
\label{CBP}

\textbf{1. General setup.}
Let $\VVV^n$, $A$, and $\AAA^n$ be the same as above.

\begin{Definition}\rm
\label{CBP1}
The set of \textit{strictly acceptable opportunities} is
$$
\AAA=\bigl\{X\in L_s^1:\exists X^1,\dots,X^N\in L_s^1
\text{ such that }\textstyle\sum_n\!X^n=X
\text{ and }X^n\in\AAA^n\;\forall n\bigr\}.
$$
\end{Definition}

\begin{Lemma}
\label{CBP2}
Let $\rho^n$ be the coherent risk measure with the
determining set~$\VVV^n$ and let $\rho$ be the convolution
of $\rho^1,\dots,\rho^N$. Then
$$
\AAA
=\conv_n\AAA^n
=\{X\in L_s^1:\rho(X)<0\}.
$$
\end{Lemma}

This statement follows from~\eqref{bo3}.

\skm
The NSAO condition, the set of risk-neutral measures,
and the upper and lower prices are defined similarly
as above.

\begin{Theorem}[FTAP]
\label{CBP3}
The NSAO condition is satisfied if and only if
$(\bigcap_n\!\VVV^n)\cap\RRR\ne\emp$.
\end{Theorem}

The proof is similar to the proof of Theorem~\ref{MBP5}.

\begin{Theorem}[Pricing]
\label{CBP4}
If $A$ is a cone and $F\in L_s^1$, then
\begin{align}
\label{cbp1}
\wl V(F)&=\sup\{\EE_\QQ F:
\QQ\in(\textstyle\bigcap_n\!\VVV^n)\cap\RRR\},\\
\label{cbp2}
\ul V(F)&=\inf\{\EE_\QQ F:
\QQ\in(\textstyle\bigcap_n\!\VVV^n)\cap\RRR\}.
\end{align}
\end{Theorem}

The proof is similar to the proof of Theorem~\ref{MBP7}.

\skm
In the maximum-based pricing, superreplicating a contingent
claim~$F$ means applying a trading strategy~$X_*$ such
that $\rho^n(X_*-F+\wl V(F))\le0$ for any~$n$.
The motivation is that the liability $X_*-F+\wl V(F)$
is riskless because it can be sold to any group at the
price~0.
For the convolution-based pricing, superreplication is
a more complicated procedure consisting of two steps:
\begin{mitemize}
\item[\bf 1.] Apply a trading strategy $X_*$ such that
$\rho(X_*-F+\wl V(F))=0$.
\item[\bf 2.] Split the liability $X_*-F+\wl V(F)$ into
contracts $Y_*^1,\dots,Y_*^N$ such that $\rho^n(Y_*^n)=0$
for any~$n$ and sell the $n$-th contract to the $n$-th
group.
\end{mitemize}
To sum up, a superreplication strategy is a collection
$X_*,Y_*^1,\dots,Y_*^N$.

\skb
\textbf{2. Theoretical example.}
Let us present the explicit solution of the pricing and
hedging problem in the case, where $\VVV^n$ is the
determining set of~$\rho_{\mu^n}$.
According to Example~\ref{BO3}~(i), $\rho=\rho_\mu$ with
$\mu(dx)=-x\Psi''(x)$, where $\Psi=\min_n\!\Psi_{\mu^n}$.
Therefore,
$$
\wl V(F)=-\inf_{X\in A}\rho_\mu(X-F).
$$

\begin{Theorem}
\label{CBP5}
{\bf(i)} Suppose that there exists
$$
X_*\in\argmin_{X\in A}\rho_\mu(X-F).
$$
Choose functions $h^1,\dots,h^N\colon[0,1]\to[0,1]$
such that $\sum_n h^n\equiv1$ and $h^n=0$ outside the
set $\{\Psi_{\mu^n}=\Psi\}$. Consider
\begin{align*}
f^n(x)&=\int_0^xh^n(f(y))dy,\quad x\in\R,\\
Y_*^n&=f^n(X_*-F+\wl V(F)),
\end{align*}
where $f$ is the distribution function of~$X_*-F+\wl V(F)$
and $\int_0^x$ is the oriented integral.
Then $X_*,Y_*^1,\dots,Y_*^N$ is a superreplication
strategy.

{\bf(ii)} If the support of each $\mu^n$ is $[0,1]$,
then $X_*$ {\rm(}if it exists{\rm)} is unique and
any superreplication strategy has the form described
above.
\end{Theorem}

\Remark
A geometric recipe for finding $X_*$ can be found
in~\cite[Subsect.~2.5]{C062}.\End

\skm
In most typical situations, the sets
$\{\Psi_{\mu^n}=\Psi\}$ are disjoint intervals
(some of them might be empty)
and $X_*-F$ has a continuous distribution.
Let $\{\Psi_{\mu^n}=\Psi\}=[a^n,b^n]$. Then
$$
f^n(x)
=\int_0^x I(f(y)\in[a^n,b^n])dy
=\int_0^x I(y\in[f^{-1}(a^n),f^{-1}(b^n)]dy.
$$
In this case the form of $Y_*^n$ (see Figure~2) is very
similar to the structure of a CDO (for the definition of
this contract, see~\cite[Sect.~9.1.3]{EFM05}).

\begin{figure}
\begin{picture}(150,50)(-64,-15)
\put(-0.3,0){\includegraphics{pricing.1}}
\put(0,0){\vector(1,0){35}}
\put(0,0){\vector(0,1){35}}
\multiput(0,30.4)(1.93,0){16}{\line(1,0){1}}
\multiput(30,0)(0,1.94){16}{\line(0,1){1}}
\multiput(12,0)(0,2.1){12}{\line(0,1){1}}
\put(-3,29){\small$1$}
\put(-1,-4){\scalebox{0.8}{$a^2$}}
\put(8,-4){\scalebox{0.8}{$b^2\!=\!a^1$}}
\put(29,-4){\scalebox{0.8}{$b^1$}}
\put(0,20){\small $\Psi_{\mu^1}$}
\put(18,32){\small $\Psi_{\mu^2}$}
\put(13,22){\small $\Psi$}
\put(8,-12){\small\textbf{Figure~1}}
\end{picture}

\begin{picture}(150,63)(-80,-35)
\put(-45.5,-20.4){\includegraphics{pricing.2}}
\put(-50,0){\vector(1,0){100}}
\put(0,-22.5){\vector(0,1){45}}
\multiput(-25,0)(0,-2){8}{\line(0,-1){1}}
\multiput(-10,0)(0,-2){5}{\line(0,-1){1}}
\multiput(15,0)(0,2){8}{\line(0,1){1}}
\multiput(25,0)(0,2){5}{\line(0,1){1}}
\put(-34.5,-6){\small $f^1$}
\put(-19.5,-6){\small $f^2$}
\put(1.5,5){\small $f^3$}
\put(16.5,5){\small $f^4$}
\put(26.5,5){\small $f^5$}
\put(-22.5,-32.5){\small\textbf{Figure~2.} The form of $f^n$}
\end{picture}
\end{figure}

\skm
\textbf{3. Empirical estimation.}
The above example provides rather an explicit solution
of the superreplication problem for the case, where each
$\VVV^n$ is the determining set of Weighted V@R.
However, as shown by Table~1, there exist many other
natural examples of valuation measures.
Here we will consider possible ways to estimate
$\wl V(F)$ empirically.
The empirical estimation of $\wl V(F)$ directly
by~\eqref{cbp1} might be problematic because it may
be hard to capture the set $(\bigcap_n\!\VVV^n)\cap\RRR$.
Therefore, instead of trying to estimate $\wl V(F)$, we
will provide an upper estimate of this value.
By~\eqref{cbp1},
$$
\wl V(F)\le\inf_{n=1,\dots,N}\sup\{\EE_\QQ F:\QQ\in\VVV^n\cap\RRR\}.
$$
The intersection $\VVV^n\cap\RRR$ is still rather an
unpleasant object. However, we can get rid of~$\RRR$
using the following equality whose proof can be found
in~\cite[Th.~2.6]{C062}:
$$
\sup\{\EE_\QQ F:\QQ\in\VVV^n\cap\RRR\}
=\inf_{X\in A}\sup\{\EE_\QQ(F-X):\QQ\in\VVV^n\}.
$$
As we are trying to find an upper estimate of $\wl V(F)$,
we need not find this infimum, but can just take the
minimum over several $X$s.

In the example below, we discuss how to find/estimate
$\sup\{\EE_\QQ X:\QQ\in\VVV\}$ (to simplify the notation,
we have replaced $F-X$ by~$X$) for the valuation
measures from Table~1.

\begin{Example}\rm
\label{CBP6}
\textbf{(i)} If $\VVV=\DDD_\mu$, then
$$
\sup\{\EE_\QQ X:\QQ\in\VVV\}=\rho_\mu(-X).
$$
A theoretical representation of $\rho_\mu$ is provided
by~\eqref{crm2}.
Its empirical estimate is provided by~\eqref{crm4}.

\textbf{(ii)} If $\VVV=\EE(\DDD_\mu\cond Y)$,
then, according to~\eqref{bo1},
$$
\sup\{\EE_\QQ X:\QQ\in\VVV\}
=\rho_\mu^f(-X;Y)
=\rho_\mu(-f(X)),
$$
where $f(y)=\EE(X\cond Y=y)$.
Thus, we should calculate $f$ and then apply the
procedures of~(i).

\textbf{(iii)} If $\VVV=\EX_{\DDD_\mu}(W)$ and $W$ has
a continuous distribution, then, according to
Example~\ref{BO1}~(i),
$$
\sup\{\EE_\QQ X:\QQ\in\VVV\}
=\EE_{\QQ_\mu(W)}X
=\EE\psi_\mu(W)X
=\rho_\mu^c(-X;W).
$$
In order to find an empirical estimate of this quantity,
take time series $(x_1,w_1),\dots,(x_T,w_T)$ for $(X,W)$.
As the unit time interval here equals the duration of
the contingent claim we are trying to price, i.e. it has
the order of several months, the ordinary time series
might not be available, and one could use the bootstrap
technique. Let $w_{(1)},\dots,w_{(T)}$ be the numbers
$w_1,\dots,w_T$ in the increasing order and let $n(i)$ be the
number such that $w_{(i)}=w_{n(i)}$. Then, according
to Example~\ref{BO1}~(ii), an empirical estimate of
$\rho_\mu^c(-X;W)$ is given by
$$
\sum_{t=1}^T x_{n(t)}\int_{(t-1)/T}^{t/T}\psi_\mu(x)dx.
$$

\textbf{(iv)} If $\VVV=\EX_{\EE(\DDD_\mu\cond Y)}(W)$,
then, according to Example~\ref{BO2} and~\eqref{bo2},
\begin{align*}
\sup\{\EE_\QQ X:\QQ\in\VVV\}
&=\sup\{\EE_\QQ f(Y):\QQ\in\VVV\}\\
&=\sup\Bigl\{\int_{\R^d}f(y)\phi(y)\wt\PP(dy):
\phi\in\wt\EX_{\wt\DDD_\mu}(g)\Bigr\}\\
&=\rho_\mu^c(f(Y);g(Y))\\
&=\rho_\mu^{fc}(X;Y;W),
\end{align*}
where $f(y)=\EE(X\cond Y=y)$ and the other notation is
introduced in Example~\ref{BO2}.
Thus, we should calculate $f$ and $g$ and then apply
the procedures of~(iii).

\textbf{(v)} If $\VVV=\conv_k\VVV^k$, then
$$
\sup\{\EE_\QQ X:\QQ\in\VVV\}
=\max_k\sup\{\EE_\QQ X:\QQ\in\VVV^k\}.
$$

\vspace{-8mm}
\hfill$\Box$
\end{Example}

\vspace{3mm}
\Remark
If $\VVV=\EX_\DDD(W)=\{\QQ\}$, where $\DDD$ is the
determining set of a risk measure~$\rho$, then
$\EE_\QQ F$ is the price of~$F$ obtained through the
coherent optimality pricing technique
of~\cite[Subsect.~3.3]{C062}, i.e. it is the number~$x$
such that
$$
\inf_{h\in\R}\rho(W+h(F-x))=\rho(W).
$$
In other words, this is the coherent reservation
price of~$F$.\End

\skm
The paper~\cite{CM061} contains more
information on the empirical estimates of
$\rho$, $\rho^f$, $\rho^c$, and $\rho^{fc}$.
In particular, it describes convenient Monte Carlo
estimation procedures for Alpha V@R and Beta V@R.

\skb
\textbf{4. Liquidity.}
The definition below was given in~\cite{C062} following
the ideas of~\cite{CGM01}.

\begin{Definition}\rm
\label{CBP7}
The \textit{upper} and \textit{lower price functions}
of a contingent claim~$F$ are defined as
\begin{align*}
\wl V(F,v)&=\inf\{x\in\R:\exists X\in A\text{ such that }
X+v(-F+x)\in\AAA\},\quad v>0,\\
\ul V(F,v)&=\sup\{x\in\R:\exists X\in A\text{ such that }
X+v(F-x)\in\AAA\},\quad v>0.
\end{align*}
From the financial point of view, $v$ means the volume
of a trade.
\end{Definition}

If $A$ is a cone, which means that there are no liquidity
effects, then $\wl V(F,v)\equiv\wl V(F)$,
$\ul V(F,v)\equiv\ul V(F)$.
However, if $A$ is bounded in a certain sense, then
the upper and lower price functions are not constant.
In view of the equality $\ul V(F,v)=-\wl V(-F,v)$,
it is sufficient to study only the properties of
$\wl V(F,v)$.

\begin{Theorem}
\label{UPL2}
Let $F\in L_s^1$.

{\bf(i)} The function $\wl V(F,\,\cdot\,)$ is
increasing and continuous.

{\bf(ii)} We have
$$
\lim_{v\da0}\wl V(F,v)
=\sup_{\QQ\in(\bigcap_n\!\VVV^n)\cap\RRR}\EE_\QQ F.
$$

{\bf(iii)} We have
$$
\lim_{v\to\infty}\wl V(F,v)
\le\sup_{\QQ\in\bigcap_n\!\VVV^n}\EE_\QQ F.
$$
If $\sup_{X\in A,\,\QQ\in\bigcap_n\!\VVV^n}|\EE_\QQ X|<\infty$,
then
$$
\lim_{v\to\infty}\wl V(F,v)
=\sup_{\QQ\in\bigcap_n\!\VVV^n}\EE_\QQ F.
$$
\end{Theorem}

This statement follows from the results
of~\cite[Subsect.~2.8]{C062}.

\section{Sensitivity Coefficients}
\label{SENS}

In this section, we adjust the above technique to measure
sensitivity coefficients.
The basic idea is as follows.
The market is valuing any contingent claim by a valuation
measure~$\QQ$. We do not know~$\QQ$ completely, but we can
indicate a set~$\VVV$, to which it belongs.
For example, the arguments of Section~\ref{CBP} say
that it should belong to $\VVV:=(\bigcap_n\!\VVV^n)\cap\RRR$.
Suppose moreover that the payoff function of a contingent
claim has the form $F=f(S,\xi)$, where $S$ is the
value of the underlying asset, with respect to which we
wish to take the sensitivity coefficient, and $\xi$ is
a random variable that does not depend on~$S$
(natural examples are given below).
Then we can say that the sensitivity coefficient is
$$
\Delta=\frac{\partial}{\partial S}\,\EE_\QQ f(S,\xi)
=\EE_\QQ\frac{\partial f}{\partial S}(S,\xi).
$$
Recalling that we do not know $\QQ$ completely, but we
know only~$\VVV$, we can define the interval for deltas
as
$$
I(\Delta):=\Bigl\{\EE_\QQ\frac{\partial f}{\partial S}\,
(S,\xi):\QQ\in\VVV\Bigr\}.
$$
The problem of finding/estimating this interval is
exactly the problem we have considered above with
$F$ replaced by $\frac{\partial f}{\partial S}$.

If the model is complete (like the
Black--Scholes--Merton one), then $\VVV=\{\QQ_0\}$ and
$\EE_{\QQ_0}f(S,\xi)=\phi(S)$ is the fair price of~$F$.
In this case
$$
I(\Delta)
=\Bigl\{\EE_{\QQ_0}\frac{\partial f}{\partial S}\,(S,\xi)\Bigr\}
=\{\phi'(S)\},
$$
which is the traditional sensitivity.
Thus, in complete models we have one fair price and
one sensitivity; in incomplete models we have an interval of
fair prices and an interval for sensitivities.
In particular, risk in incomplete models cannot
be completely eliminated by the delta hedging.

\begin{Example}\rm
\label{SENS1}
Let $F$ be the standard call option, i.e.
$F=e^{-rT}(S_T-K)^+$, where $S_T$ is the price of some
asset at the expiration date~$T$ and $r$ is the risk-free
rate. It is natural to model $S_T$ as $S\xi$, where $S$
is the current price of the asset and $\xi$ is a random
variable whose distribution does not depend on~$S$.
Then $f(S,\xi)=e^{-rT}(S\xi-K)^+$ and
$$
\frac{\partial f}{\partial S}\,(S,\xi)
=e^{-rT}\xi I(\xi\ge K/S).
$$

\vspace{-8mm}
\hfill$\Box$
\end{Example}

\vspace{3mm}
\begin{Example}\rm
\label{SENS2}
Let $F$ be the option on a bond providing the amounts
$c_1,\dots,c_N$ at future dates $T_1<\dots<T_N$.
The expiration date of the option is $T<T_1$.
Suppose that the yield curve evolves as
$r(t,T)=r_t+\phi(T-t)$, where $r_t$ is a random short
rate and $\phi$ is a fixed shape (for example, this is
the case for the Vasicek model). Then the discounted
payoff of the option is
$$
F=e^{-Tr(0,T)}\Bigl(\sum_{n=1}^N c_nf_n(r_t)-K\Bigr)^+,
$$
where
$$
f_n(r)=\exp\{-(T_n-T)r-(T_n-T)\phi(T_n-T)\}.
$$
If we again model $r_t$ as $r\xi$, where $r$ is the current
short rate and the distribution of $\xi$ does not depend
on~$r$, then $F=f(r,\xi)$ and
$$
\frac{\partial f}{\partial r}\,(r,\xi)
=e^{-Tr(0,T)}\sum_{n=1}^Nc_nf'_n(r\xi)I\Bigl(\sum_{n=1}^N
c_n f_n(r\xi)\ge K\Bigr).
$$

\vspace{-10mm}
\hfill$\Box$
\end{Example}

\section{Summary and Conclusion}
\label{SC}

In this paper, we propose a technique for pricing and
hedging in incomplete markets that yields finer price
intervals than the technique proposed in~\cite{CGM01}.
The basic idea is: to hedge a liability means to employ
a trading strategy and then split the resulting liability
into several contracts that can be sold to different
groups at the price~0.
The corresponding interval of fair prices is
\begin{equation}
\label{sc1}
I(F)=\{\EE_\QQ F:\QQ\in(\textstyle\bigcap_n\!\VVV^n)\cap\RRR\}.
\end{equation}
Here $\VVV^n$ is the set of valuation measures used by
the $n$-th representative agent (for example, these
agents are large companies) or a group of agents.
A pleasant feature of this technique, which is not shared
by the technique of~\cite{CGM01}, is: the more groups
are taken into account, the smaller are the fair price
intervals.

The intersection $(\bigcap_n\!\VVV^n)\cap\RRR$ is nonempty
provided that there exists no trade involving
all the groups, after which each group is better off
(Fundamental Theorem of Asset Pricing).

For the estimation purposes, we can
replace~\eqref{sc1} by a wider interval
$$
\bigcap_{n=1}^N\{\EE_\QQ F:\QQ\in\VVV^n\cap\RRR\}.
$$
For typical choices of valuation measures
(see Table~1), the latter interval admits simple
empirical estimation procedures (see Example~\ref{CBP6}).

The interval~\eqref{sc1} can be made very small if we
know the portfolios of the representative agents.
If the $n$-th agent is using the expected utility to
assess the quality of his/her position, then his/her
valuation measure is $cU'(W_1)\PP$,
where $U$ is the utility function and $W_1$ is the
terminal wealth.
If the agent is employing the coherent risk/utility
measurement, then a substitute for the above measure
is the extreme measure.
If, however, we have no information on the structure of
the agents' portfolios, then a natural choice of the
set of valuation measures is the determining set of
some coherent risk measure.
Both for the risk measurement purposes and for the
estimation of fair price intervals, it is convenient
to use not pure coherent risk measures, but rather their
factor versions introduced in~\cite{CM061}.

Furthermore, we define the interval of sensitivities
of a contingent claim $F=f(S,\xi)$ with respect to the
price~$S$ of an underlying asset as
$$
I(\Delta)=\Bigl\{\EE_\QQ\frac{\partial f}{\partial S}\,
f(S,\xi):\QQ\in(\textstyle\bigcap_n\!\VVV^n)\cap\RRR\Bigr\}.
$$
If we manage to make $(\bigcap_n\!\VVV^n)\cap\RRR$ small,
then both the fair price intervals and the sensitivity
intervals are small; in particular, the risk of holding~$F$
can be successfully eliminated by the delta hedging.

\clearpage

\end{document}